\input amstex 
\documentstyle{amsppt}
\input bull-ppt
\keyedby{bull335/lbd}

\topmatter
\cvol{28}
\cvolyear{1993}
\cmonth{January}
\cyear{1993}
\cvolno{1}
\cpgs{84-89}
\title Linkless embeddings of graphs in $3$-space\endtitle
\author Neil Robertson, P. D. Seymour, and Robin 
Thomas\endauthor
\shortauthor{N. Robertson, P. D. Seymour, and R. Thomas}
\shorttitle{Linkless embeddings of graphs in 3-space}
\address Department of Mathematics, Ohio State University, 
231 West 18th
Avenue,
Columbus, Ohio 43210\endaddress
\ml robertso\@function.mps.ohio-state.edu\endml
\address Bellcore, 445 South Street, Morristown, New 
Jersey 07962\endaddress
\ml pds\@bellcore.com\endml
\address School of Mathematics, Georgia Institute of 
Technology, Atlanta,
Georgia 30332\endaddress
\ml thomas\@math.gatech.edu\endml
\date January 14, 1992 and, in revised form, May 12, 
1992\enddate
\thanks The first author's research was performed under a 
consulting agreement
with Bellcore; he was supported by NSF under Grant No. 
DMS-8903132 and by ONR
under Grant No. N00014-911-J-1905. The third author was 
supported in part by
NSF under Grant No. DMS-9103480, and in part by DIMACS 
Center, Rutgers
University, New Brunswick, New Jersey 08903\endthanks
\abstract We announce results about flat (linkless) 
embeddings of graphs in
3-space. A piecewise-linear embedding of a graph in 
3-space is called {\it flat\/}
if every circuit of the graph bounds a 
disk disjoint from the rest of
the graph. We have shown:
\par
(i) An embedding is flat if and only if the fundamental 
group of the complement
in 3-space of the embedding of every subgraph is free.
\par
(ii) If two flat embeddings of the same graph are not 
ambient isotopic, then
they differ on a subdivision of $K_5$ or $K_{3,3}$.
\par
(iii) Any flat embedding of a graph can be transformed to 
any other flat
embedding of the same graph by ``3-switches'', an analog 
of 2-switches from the
theory of planar embeddings. In particular, any two flat 
embeddings of a
4-connected graph are either ambient isotopic, or one is 
ambient isotopic to a
mirror image of the other.
\par
(iv) A graph has a flat embedding if and only if it has no 
minor isomorphic to
one of seven specified graphs. These are the graphs that 
can be obtained from
$K_6$ by means of $Y\Delta$- and $\Delta Y$-exchanges.\endabstract
\endtopmatter

\document
\heading 1. Introduction\endheading
\par
All spatial embeddings are assumed to be piecewise linear. 
If $C,C'$ are
disjoint simple closed curves in $S^3$, then their {\it 
linking number\/},
$\roman{lk}(C,C')$, is the number of times (mod 2) that 
$C$ crosses over $C'$
in a regular projection of  $C\cup C'$. In this paper 
graphs are finite,
undirected, and may have loops and multiple edges. Every 
graph is regarded as a
topological space in the obvious way. We say that an 
embedding of a graph $G$
in $S^3$ is {\it linkless\/} if every two disjoint 
circuits of $G$ have zero
linking number. The following is a result of Sachs [13, 
14] and Conway and
Gordon \cite3.
\thm\nofrills{\rm{(1.1)}}\quad The graph $K_6$ {\rm(}the 
complete graph on six vertices{\rm)}
has no linkless embedding.
\ethm
\demo{Proof} Let $\phi$ be an embedding of $K_6$ into 
$S^3$. By studying the
effect of a crossing change in a regular projection, it is 
easy to see that the
mod 2 sum $\sum\roman{lk}
(\phi(C_1),\phi(C_2))$, where the sum is taken over all 
unordered pairs of
disjoint circuits $C_1,C_2$ of $K_6$, is an invariant 
independent of the
embedding. By checking an arbitrary embedding we can 
establish that this invariant
equals 1.\qed\enddemo
\par
Let $G$ be a graph and let $v$ be a vertex of $G$ of 
valency 3 with distinct
neighbors. Let $H$ be obtained from $G$ by deleting $v$ 
and adding an edge
between every pair of neighbors of $v$. We say that $H$ is 
obtained from $G$ by
a $Y\Delta$-{\it exchange\/} and that $G$ is obtained from 
$H$ by a $\Delta
Y$-{\it exchange\/}. The {\it Petersen family\/} is the 
set of all graphs that
can be obtained from $K_6$  by means of $Y\Delta$- and 
$\Delta Y$-exchanges.
There are exactly seven such graphs, one of which is the 
Petersen graph. Pictures
of these graphs can be found in [13--15]. Sachs [13, 14] 
has in fact shown that
no member of the Petersen family has a linkless embedding 
[the argument is
similar to the proof of (1.1)] and raised the problem of 
characterizing
linklessly embeddable graphs. A graph is a {\it minor\/} 
of another if the
first can be obtained from a subgraph of the second by 
contracting edges. It is
easy to see that the property of having a linkless 
embedding is preserved under
taking minors,  and that led Sachs to conjecture that a 
graph is linklessly
embeddable if and only if it has no minor in the Petersen 
family. We have shown
that this is true. Moreover, let us say
that an embedding $\phi$ of a graph $G$ in $S^3$ is {\it 
flat\/} if for every
circuit $C$ of $G$ there exists an open disk in $S^3$ 
disjoint from $\phi(G)$
whose boundary is $\phi(C)$. Clearly every flat embedding 
is linkless, but the
converse need not hold. However, B\"ohme \cite1 and Saran 
\cite{15} conjectured
that a graph has a linkless embedding if and only if it 
has a flat one. This is
also true, for we have shown the following.
\thm\nofrills{\rm{(1.2)}}\quad For a graph $G$, the 
following are equivalent{\rm:}
\roster
\item "(i)" $G$ has a flat embedding,
\item "(ii)" $G$ has a linkless embedding,
\item "(iii)" $G$ has no minor in the Petersen family.
\endroster
\ethm
\par
There have been a number of other attempts [8, 15, 2] at 
proving
$\roman{(iii)}\Rightarrow\roman{(i)}$ and
$\roman{(iii)}\Rightarrow\roman{(ii)}$. However, none of 
them is correct.
\par
For the proof of (1.2) we need the following two theorems, 
which may be of
independent interest.
\thm\nofrills{\rm{(1.3)}}\quad Let $\phi$ be an embedding 
of a graph $G$ in $S^3$. Then
$\phi$ is flat if and only if for every subgraph $G'$ of 
$G$, the fundamental
group of $S^3-\phi(G')$ is free.
\ethm
\par
Let $\phi_1,\phi_2$ be two embeddings of a graph $G$ in 
$S^3$. We say that
$\phi_1,\phi_2$ are {\it ambient isotopic\/} if there 
exists an orientation
preserving homeomorphism $h$ of $S^3$ onto $S^3$ such that 
$\phi_1=h\phi_2$.
(We remark that by a result of Fisher \cite4 $h$ can be 
realized by an ambient
isotopy.) If $\phi$ is an embedding of a graph $G$ in 
$S^3$ we denote by
$-\phi$ the embedding of $G$ obtained by composing $\phi$ 
with the antipodal
map.
\thm\nofrills{\rm{(1.4)}}\quad Let $G$ be a $4$-connected 
graph and let $\phi_1,\phi_2$ be
two flat embeddings of $G$. Then $\phi_1$ is ambient 
isotopic to either
$\phi_2$ or $-\phi_2$.
\ethm
\heading 2. The fundamental group\endheading
\par
A basic tool for working with flat embeddings is the 
following lemma of B\"ohme
\cite1 (see also \cite{15}).
\thm\nofrills{\rm{(2.1)}}\quad Let $\phi$ be a flat 
embedding of a graph $G$ into $S^3$, and
let $C_1,C_2,\dotsc,C_n$ be a family of circuits of $G$ 
such that for every
$i\not=j$, the intersection of $C_i$ and $C_j$ is either 
connected or null. Then
there exist pairwise disjoint open disks 
$D_1,D_2,\dotsc,D_n$, disjoint from
$\phi(G)$ and such that $\phi(C_i)$ is the boundary of 
$D_i$ for
$i=1,2,\dotsc,n$.
\ethm
\par
We illustrate the use of (2.1) with the following, which 
is a special case of a
theorem of Wu \cite{18}. An embedding $\phi$ of a graph 
$G$ in $S^3$ is {\it
spherical\/} if there exists a surface $\Sigma\subseteq 
S^3$ homeomorphic to
$S^2$ such that $\phi(G)\subseteq\Sigma$. Clearly if 
$\phi$ is spherical then
$G$ is planar.
\thm\nofrills{\rm{(2.2)}}\quad Let $\phi$ be an embedding 
of a planar graph $G$ in $S^3$.
Then $\phi$ is flat if and only if it is spherical.
\ethm
\demo{Proof} Clearly if $\phi$ is spherical then it is 
flat. We prove the
converse only for the case when $G$ is 3-connected. Let 
$C_1,C_2,\dotsc,C_n$
be the collection of face-boundaries in some planar 
embedding of $G$. These
circuits satisfy the hypothesis of (2.1). Let 
$D_1,D_2,\dotsc,D_n$ be the disks
as in (2.1); then $\phi(G)\cup D_1\cup D_2\cup\dotsb\cup 
D_n$ is the desired
sphere.\qed\enddemo
\par
The following is a result of Scharlemann and Thompson 
\cite{16}.
\thm\nofrills{\rm{(2.3)}}\quad Let $\phi$ be an embedding 
of a graph $G$ in $S^3$. Then
$\phi$ is spherical if and only if
\roster
\item "{(i)}" $G$ is planar, and
\item "{(ii)}" for every subgraph $G'$ of $G$, the 
fundamental group of
$S^3-\phi(G')$ is free.
\endroster
\ethm
\par
We see that by (2.2), (1.3) is a generalization of (2.3). 
In fact, we prove
(1.3) by reducing it to planar graphs and then applying 
(2.3). Let us prove
the ``only if'' part of (1.3). Let $G'$ be a subgraph of 
$G$ such that
$\pi_1(S^3-\phi(G'))$ is not free. Choose a maximal forest 
$F$ of $G'$ and let
$G''$ be obtained from $G'$ by contracting all edges of 
$F$, and let $\phi''$ be
the induced embedding of $G''$. Then
$\pi_1(S^3-\phi''(G''))=\pi_1(S^3-\phi(G'))$ is not free, 
but $G''$ is planar,
and so $\phi''$ is not flat by (2.2) and (2.3). Hence 
$\phi$ is not flat, as
desired.
\par
Let $G$ be a graph, and let $e$ be an edge of $G$. We 
denote by $G\backslash
e(G/e)$ the graph obtained from $G$ by deleting 
(contracting) $e$. If $\phi$ 
is an embedding of $G$ in $S^3$, then it induces 
embeddings of $G\backslash e$
and (up to ambient isotopy) of $G/e$ in the obvious way. 
We denote these
embeddings by $\phi\backslash e$ and $\phi/e$, respectively.
\thm\nofrills{\rm{(2.4)}}\quad Let $\phi$ be an embedding 
of a graph $G$ into $S^3,$ and let
$e$ be a nonloop edge of $G$. If both $\phi\backslash e$ 
and $\phi/e$ are flat,
then $\phi$ is flat.
\ethm
\demo{Proof} Suppose that $\phi$ is not flat. By (1.3) 
there exists a subgraph
$G'$ of $G$ such that $\pi_1(S^3-\phi(G'))$ is not free. 
If $e\not\in E(G')$
then $\phi\backslash e$ is not flat by (1.3). If $e\in 
E(G')$ then $\phi/e$ is
not flat by (1.3), because 
$\pi_1(S^3-(\phi/e)(G'/e))=\pi_1(S^3-\phi(G'))$ is
not free.\qed\enddemo
\par
We say that a graph $G$ is a {\it coforest\/} if every 
edge of $G$ is a loop.
The following follows immediately from (2.4).
\thm\nofrills{\rm{(2.5)}}\quad Let $\phi$ be an embedding 
of a graph $G$ in $S^3$. Then
$\phi$ is flat if and only if the induced embedding of 
every coforest minor of
$G$ is flat.
\ethm
\heading 3. Uniqueness\endheading
\par
A graph $H$ is a {\it subdivision\/} of a graph $G$ if $H$ 
can be obtained from
$G$ by replacing edges by pairwise internally-disjoint 
paths. We recall that
Kuratowski's theorem \cite6 states that a graph is planar 
if and only if it
contains no subgraph isomorphic to a subdivision of $K_5$ 
or $K_{3,3}$. It
follows from a theorem of Mason \cite7 and (2.2) that any 
two flat embeddings
of a planar graph are ambient isotopic. On the other hand 
we have the
following.
\thm\nofrills{\rm{(3.1)}}\quad The graphs $K_5$ and 
$K_{3,3}$ have exactly two nonambient
isotopic flat embeddings.
\ethm
\demo{Sketch of proof} Let $G$ be $K_{3,3}$ or $K_5$, let 
$e$ be an edge of
$G$, and let $H$ be $G\backslash e$. Notice that $H$ is 
planar. From (2.1) it
follows that if $\phi$ is a flat embedding of $G$, then 
there is an embedded
2-sphere $\Sigma\subseteq S^3$ with 
$\phi(G)\cap\Sigma=\phi(H)$.
If $\phi_1$ and $\phi_2$ are flat embeddings of $G$, we 
may assume
(by replacing $\phi_2$ by
an ambient isotopic embedding) that this 2-sphere $\Sigma$ 
is the same for both
$\phi_1$ and $\phi_2$. Now $\phi_1$ is ambient isotopic to 
$\phi_2$ if and only
if  $\phi_1(e)$ and $\phi_2(e)$ belong to the same 
component of
$S^3-\Sigma$.\qed\enddemo
\par
As a curiosity we deduce that a graph has a unique flat 
embedding if and only
if it is planar.
\par
We need the following three lemmas. We denote by $f|X$ the 
restriction of a
mapping $f$ to a set $X$.
\thm\nofrills{\rm{(3.2)}}\quad Let $\phi_1,\phi_2$ be two 
flat embeddings of a graph $G$ that
are not ambient isotopic. Then there exists a subgraph $H$ 
of $G$ isomorphic to
a subdivision of $K_5$ or $K_{3,3}$ for which $\phi_1|H$ 
and $\phi_2|H$ are not
ambient isotopic.
\ethm
\par
We denote the vertex-set and edge-set of a graph $G$ by 
$V(G)$ and $E(G)$
respectively. Let $G$ be a graph and let $H_1,H_2$ be 
subgraphs of $G$
isomorphic to subdivisions of $K_5$ or $K_{3,3}$. We say 
that $H_1$ and $H_2$
are 1-{\it adjacent\/} if there exist $i\in\{1,2\}$ and a 
path $P$ in $G$ such
that $P$ has only its endvertices in common with $H_i$ and 
such that $H_{3-i}$
is a subgraph of the graph obtained from $H_i$ by adding 
$P$. We say that
$H_1$ and $H_2$ are 2-{\it adjacent\/} if there are seven 
vertices
$u_1,u_2,\dotsc,u_7$ of $G$ and thirteen paths $L_{ij}$ of 
$G$ $(1\le i\le 4$
and $5\le j\le 7$, or $i=3$ and $j=4)$, such that
\roster
\item "(i)" each path $L_{ij}$ has ends $u_i,u_j$,
\item "(ii)" the paths $L_{ij}$ are mutually 
vertex-disjoint except for their
ends,
\item "(iii)" $H_1$ is the union of $L_{ij}$ for $i=2,3,4$ 
and $j=5,6,7$, and
\item "(iv)" $H_2$ is the union of $L_{ij}$ for $i=1,3,4$ 
and $j=5,6,7$.
\endroster
(Notice that if $H_1$ and $H_2$ are 2-adjacent then they 
are both isomorphic to
subdivisions of $K_{3,3}$ and that $L_{34}$ is used in 
neither $H_1$ nor
$H_2.)$ We denote by $\scr K(G)$ the simple graph with 
vertex-set all subgraphs
of $G$ isomorphic to subdivisions of $K_5$ or $K_{3,3}$ in 
which two distinct
vertices are adjacent if they are either 1-adjacent or 
2-adjacent. The
following is easy to see, using (3.1).
\thm\nofrills{\rm{(3.3)}}\quad Let $\phi_1,\phi_2$ be two 
flat embeddings of a graph $G$, and
let $H,H'$ be two adjacent vertices of $\scr K(G)$. If 
$\phi_1|H$ is ambient
isotopic to $\phi_2|H$, then $\phi_1|H'$ is ambient 
isotopic to $\phi_2|H'$.
\ethm
\par
The third lemma is purely graph-theoretic.
\thm\nofrills{\rm{(3.4)}}\quad If $G$ is a $4$-connected 
graph, then $\scr K(G)$ is
connected.
\ethm
\par
We prove (3.4) in \cite{10} by proving a stronger result, 
a necessary and
sufficient condition for $H,H'\in V(\scr K(G))$ to belong 
to the same component
of $\scr K(G)$ in an arbitrary graph $G$. The advantage of 
this approach is
that it permits an inductive proof using the techniques of 
deleting and
contracting edges.
\demo{Proof of {\rm(1.4)}} If $G$ is planar then $\phi_1$ 
is ambient isotopic
to $\phi_2$ by Mason's theorem. Otherwise there exists, by 
Kuratowski's
theorem, a subgraph $H$ of $G$ isomorphic to a subdivision 
of $K_5$ or
$K_{3,3}$. By replacing $\phi_2$ by $-\phi_2$ we may 
assume by (3.1) that
$\phi_1|H$ is ambient isotopic to $\phi_2|H$. From (3.3) 
and (3.4) we deduce
that $\phi_1|H'$ is ambient isotopic to $\phi_2|H'$ for 
every $H'\in V(\scr
K(G))$. By (3.2) $\phi_1$ and $\phi_2$ are ambient 
isotopic, as desired.\qed\enddemo
\par
We now state a generalization of (1.4). Let $\phi$ be a 
flat embedding of a
graph $G$, and let $\Sigma\subseteq S^3$ be a surface 
homeomorphic to $S^2$
meeting $\phi(G)$ in a set $A$ containing at most three 
points. In one of the
open balls into which $\Sigma$ divides $S^3$, say $B$, 
choose an open disk $D$
with boundary a simple closed curve $\partial D$ such that 
$A\subseteq\partial
D\subseteq\Sigma$. Let $\phi'$ be an embedding obtained 
from $\phi$ by taking a
reflection of $\phi$ through $D$ in $B$ and leaving $\phi$ 
unchanged in
$\Sigma-B$. We say that $\phi'$ is obtained from $\phi$ by 
a 3-{\it switch\/}.
The following analog of a theorem of Whitney \cite{17} 
generalizes (1.4).
\thm\nofrills{\rm{(3.5)}}\quad Let $\phi_1,\phi_2$ be two 
flat embeddings of a graph $G$ in
$S^3$. Then $\phi_2$ can be obtained from $\phi_1$ by a 
series of $3$-switches.
\ethm
\heading 4. Main theorem\endheading
\par
The difficult part of (1.2) is to show that (iii) implies 
(i). Let us just very
briefly sketch the main idea of the proof. Suppose that 
$G$ is a minor-minimal
graph with  no flat embedding. We first show that a 
$Y\Delta$-exchange
preserves the property of having a flat embedding; thus we 
may assume that $G$
has no triangles (and indeed has some further properties 
that we shall not
specify here). It can be shown that $G$ satisfies a 
certain weaker form of
5-connectivity. Suppose that there are two edges $e,f$ of 
$G$ so that
$G\backslash e/f$ and $G/e/f$ are ``Kuratowski 
4-connected''. (Kuratowski
4-connectivity is a slight weakening of 4-connectivity for 
which (1.4) still
remains true.) Since $G$ is minor-minimal with no flat 
embedding, there are
flat embeddings $\phi_1,\phi_2,\phi_3$ of $G\backslash 
e,G/e,G/f$,
respectively. Since $G\backslash e/f$ and $G/e/f$ are both 
Kuratowski
4-connected, we can assume (by replacing $\phi_1$ or 
$\phi_2$ or both by its
mirror image) that $\phi_1/f$ is ambient isotopic to 
$\phi_3\backslash e$ and
that $\phi_2/f$ is ambient isotopic to $\phi_3/e$. Now it 
can be argued (the
details are quite complicated, see \cite{12}) that the 
uncontraction of $f$ in
$\phi_1/f\simeq\phi_3\backslash e$ is the same as in 
$\phi_2/f\simeq\phi_3/e$.
Let $\phi$ be obtained from $\phi_3$ by doing this 
uncontraction; then
$\phi\backslash e$ is ambient isotopic to $\phi_1$ and 
$\phi/e$ is ambient
isotopic to $\phi_2$. Since both these embeddings are 
flat, $\phi$ is flat by
(2.4), a contradiction. Thus no two such edges $e,f$ 
exist. But now a purely
graph-theoretic argument \cite{11} (using the nonexistence 
of such edges $e,f$,
the high connectivity of $G$, and that the graph obtained 
from $G$ by deleting
$v$ is nonplanar for every vertex $v$ of $G)$ implies $G$ 
has a minor in the
Petersen family.
\par
Finally we would like to mention some algorithmic aspects 
of flat embeddings. In
\cite{16} Scharlemann and Thompson describe an algorithm 
to test if a given
embedding is spherical. Using their algorithm, (2.2), and 
(2.5), we can test if
a given embedding is flat, by testing the flatness of all 
coforest minors. At
the moment there is no known {\it polynomial-time\/} 
algorithm to test if an
embedding of a given coforest is flat, because it includes 
testing if a knot is
trivial. On the other hand, we can test in time 
$O(|V(G)|^3)$ if a given graph
$G$ has a flat embedding. This is done by testing the 
absence of minors
isomorphic to members of the Petersen family, using the 
algorithm \cite9 of the
first two authors.

\enddocument